\documentclass[12pt]{article}

\usepackage{amsmath,amssymb}
\usepackage{fancyhdr,a4wide}
\usepackage{textcomp}
\usepackage{amssymb}
\usepackage{bbm}

\usepackage{hyperref}

\newcommand{\M}{\mathcal{M}}

\newcommand{\ol}{\overline}
\newcommand{\wt}{\widetilde}
\newcommand{\wh}{\widehat}

\newcommand{\defeq}{\stackrel{\mathrm{def.}}{=}}
\renewcommand{\epsilon}{\varepsilon}

\newcommand{\R}{\mathbb{R}}
\newcommand{\Rd}{\R^d}
\newcommand{\EXP}{{\mathbb{E}}}
\newcommand{\PROB}{\mathbb{P}}

\newcommand{\var}{\mathrm{Var}}
\newcommand{\Tr}{\mathrm{Tr}}
\newcommand{\lambdamax}{\lambda_{\text{max}}}

\newtheorem{theorem}{Theorem}

\newtheorem{proposition}{Proposition}

\begin{document}


\title{On the estimation of the mean of a random vector
}
\author{Emilien Joly 
\thanks{
Supported by the  French Agence  Nationale  de la
Recherche (ANR), under grant ANR-13-BS01-0005 (project SPADRO).
}\\
Université Paris Ouest \\
Nanterre, France; \\
emilien.joly@u-paris10.fr
 \and 
G\'abor Lugosi
\thanks{
Supported by 
the Spanish Ministry of Economy and Competitiveness,
Grant MTM2015-67304-P and FEDER, EU.
}
 \\
ICREA and Department of Economics,\\
  Pompeu Fabra University, \\
Barcelona, Spain; \\ 
gabor.lugosi@upf.edu
\and 
Roberto Imbuzeiro Oliveira
\thanks{Support from CNPq, Brazil
via {\em Ci\^{e}ncia sem Fronteiras} grant \# 401572/2014-5.
Supported by a {\em Bolsa de Produtividade em Pesquisa} from CNPq, Brazil.
Supported by FAPESP Center for Neuromathematics (grant\# 2013/ 07699-0 , FAPESP - S. Paulo Research Foundation).} \\
IMPA, Rio de Janeiro, RJ, \\
Brazil;\\ 
rimfo@impa.br
}

\maketitle

\begin{abstract}
We study the problem of estimating the mean of a multivariate
distribution based on independent samples. The main result is the 
proof of existence of an estimator with a non-asymptotic 
sub-Gaussian performance 
for all distributions satisfying some mild moment assumptions. 
\end{abstract}



\section{Introduction}

Let $X$ be a random vector taking values in $\R^d$. We
assume throughout the paper that the mean vector $\mu = \EXP X$
and covariance matrix $\Sigma= (X-\mu) (X-\mu)^T$ exist.
Given $n$ independent, identically distributed samples
$X_1,\ldots,X_n$ drawn from the distribution of $X$, one 
wishes to estimate the mean vector.

A natural and popular choice is the sample mean $(1/n)\sum_{i=1}^n X_i$
that is known to have a near-optimal behavior whenever the
distribution is sufficiently light tailed. However, whenever
heavy tails are a concern, the sample mean is to be avoided 
as it may have a sub-optimal performance. While the one-dimensional case 
(i.e., $d=1$) is quite well understood (see \cite{Cat10}, \cite{DeLeLuOl16}),
various aspects of the multidimensional problem are still to be 
revealed. This paper aims at contributing to the understanding of the
multi-dimensional case.

Before stating the main results, we briefly survey properties of
some mean estimators of real-valued random variables. Some of 
these techniques serve as basic building blocks for the estimators 
we propose for the vector-valued case. 

\subsection{Estimating the mean of a real-valued random variable}
\label{sec:1d}

When $d=1$, the simplest and most popular mean estimator is the sample mean
$\ol{\mu}_n= (1/n)\sum_{i=1}^n X_i$. The sample mean is unbiased and
the central limit theorem guarantees an asymptotically Gaussian distribution.
However, unless the distribution of $X$ has a light (e.g.,
sub-Gaussian) tail, there are no non-asymptotic sub-Gaussian performance guarantees
for $\ol{\mu}_n$. We refer the reader to Catoni \cite{Cat10} for details.
However, perhaps surprisingly, there exist estimators of $\mu$ with
much better concentration properties, see Catoni \cite{Cat10} and
Devroye, Lerasle, Lugosi, and Oliveira \cite{DeLeLuOl16}.

A conceptually simple and quite powerful estimator is the so-called
\emph{median-of-means} estimator that has been proposed, in different 
forms, in various papers, see 
Nemirovsky and Yudin \cite{NeYu83},
Hsu \cite{HsuBlog},
Jerrum, Valiant, and Vazirani \cite{JeVaVa86},
Alon, Matias, and Szegedy \cite{AMS02}. 
The median-of-means estimator is defined as follows.
Given a
  positive integer $b$ and $x_1,\ldots,x_b\in\R$, let $q_{1/2}$
  denote the median of these numbers, that is, 
$$q_{1/2}(x_1,\ldots,x_b) = x_i,\mbox{ where }\# \{k\in [b]\,:\,x_k\leq x_i\}\geq \frac{b}{2}\mbox{ and } \# \{k\in [b]\,:\,x_k\geq x_i\}\geq \frac{b}{2}.$$
(If several $i$ fit the above description, we take the smallest one.)

For any fixed $\delta\in [e^{1-n/2},1)$, first choose $b=\left\lceil \ln(1/\delta)\right\rceil$ and note that $b\leq n/2$ holds.
Next, partition $[n]=\{1,\dots,n\}$ into $b$ blocks $B_1,\ldots,B_b$, each of
size 
$|B_i|\geq \lfloor n/b\rfloor\geq 2$.
Given $X_1,\ldots,X_n$, we
compute the sample mean in each block
$$Y_i=\frac{1}{|B_i|}\sum_{j\in B_i}X_j$$
and define the median-of-means estimator by 
$\wh{\mu}_n^{(\delta)}= q_{1/2}(Y_1,\ldots,Y_B).$ 
One can show (see, e.g., Hsu \cite{HsuBlog}) that for any $n\geq 4$,
\begin{equation}
\label{eq:mom}
\PROB\left\{|\wh{\mu}_n^{(\delta)} - \mu|>2e\sqrt{2 \var(X)}
    \sqrt{\frac{(1+\ln(1/\delta))}{n}} \right\} \leq \delta~,
\end{equation}
where $\var(X)$ denotes the variance of $X$.

Note that the median-of-means estimator $\wh{\mu}_n^{(\delta)}$ does
not require any knowledge of the variance of $X$. However, it
depends on the desired confidence level $\delta$ and the partition $B_1,\dots,B_b$. Any partition satisfying $\forall i,\ |B_i|\geq \lfloor n/b\rfloor$ is valid in order to get \eqref{eq:mom}. Hence, we do not keep the dependence on the partition $B_1,\dots,B_b$ in the notation $\wh{\mu}_n^{(\delta)}$.
Devroye, Lerasle, Lugosi, and Oliveira \cite{DeLeLuOl16} introduce estimators that work for a large range of
confidence levels under some mild assumptions. Catoni \cite{Cat10}
introduces estimators of quite different flavor and gets a
non-asymptotic result of the same form as \eqref{eq:mom}. Bubeck,
Cesa-Bianchi and Lugosi \cite{BuCeLu13} apply these estimators 
in the context of bandit problems.

\subsection{Estimating the mean of random vectors}

Consider now the multi-dimensional case when $d>1$. The sample mean
$\ol{\mu}_n= (1/n)\sum_{i=1}^n X_i$ is still an obvious choice for
estimating the mean vector $\mu$.

If $X$ has a multivariate normal distribution with mean vector $\mu$
and
covariance matrix $\Sigma$, then $\ol{\mu}_n$ is also multivariate
normal with mean $\mu$ and covariance matrix $(1/n)\Sigma$ and
therefore, for $\delta \in (0,1)$, with probability at least $1-\delta$,
\begin{equation}
\label{eq:subgauss}
   \|\ol{\mu}_n-\mu\| \le \sqrt{\frac{\Tr(\Sigma)}{n}} 
    + \sqrt{\frac{2\lambdamax \log(1/\delta)}{n}}~,
\end{equation}
where $\Tr(\Sigma)$ and $\lambdamax$ denote the trace and largest
eigenvalue of the covariance matrix, respectively
(Hanson and Wright \cite{HaWr71}). For non-Gaussian and possibly
heavy-tailed distributions, one cannot expect such a sub-Gaussian 
behavior of the sample mean. The main goal of this paper is to
investigate under what conditions it is possible to define mean
estimators that reproduce a (non-asymptotic) sub-Gaussian performance
similar to (\ref{eq:subgauss}).

Lerasle and Oliveira \cite{LerasleOliveira_Robust},
Hsu and Sabato \cite{HsSa16}, and Minsker \cite{Min15} extend the
median-of-means estimator to more general spaces. In particular,
Minsker's results imply that for each $\delta\in (0,1)$ 
there exists a mean estimator $\wt{\mu}_n^{(\delta)}$ and a universal
constant $C$ such that, 
with probability at least $1-\delta$,
\begin{equation}
\label{eq:minsker}
 \|\wt{\mu}_n^{(\delta)}-\mu\| \le C\sqrt{\frac{\Tr(\Sigma) \log(1/\delta)}{n}}~.
\end{equation}
While this bound is quite remarkable--note that no assumption other
than the existence of the covariance matrix is made--, it does not
quite achieve a sub-Gaussian performance bound that resembles
\eqref{eq:subgauss}. 
An instructive example is when all eigenvalues are identical and equal to $\lambdamax$.
If the dimension $d$ is large, \eqref{eq:subgauss} is of the order of
$\sqrt{(\lambdamax/n)(d+\log (\delta^{-1}))}$ while \eqref{eq:minsker}
gives the order $\sqrt{(\lambdamax/n)(d\log (\delta^{-1}))}$.
The main result of this paper is the construction of a mean estimator that, under some mild moment
assumptions, achieves a sub-Gaussian performance bound in the sense of \eqref{eq:subgauss}. More
precisely, we prove the following.

\begin{theorem}
\label {thm:main}
For all $\delta \in (0,1)$ there exists a mean
estimator $\wh{\mu}_n^{(\delta)}$ and a universal constant $C$ such that
if $X_1,\ldots,X_n$ are i.i.d.\ random vectors in $\Rd$ with mean
$\mu\in \Rd$ and covariance matrix $\Sigma$ such that there
exists a constant $K>0$ such that, for all $v\in \Rd$ with $\|v\|=1$,
\[
   \EXP\left[ \left((X-\mu)^Tv\right)^4\right] \le K (v^T\Sigma v)^2~,
\]
then for all $n\ge CK \log d \left(d + \log (1/\delta)\right)$,
\[
   \|\wh{\mu}_n^{(\delta)}-\mu\| 
\le C  \left(\sqrt{\frac{\Tr(\Sigma)}{n}}    
 + \sqrt{\frac{\lambdamax \log(\log d/\delta)}{n}} \right)~.
\]
\end{theorem}

The theorem guarantees the existence of a mean estimator whose
performance matches the sub-Gaussian bound \eqref{eq:subgauss}, 
up to the additional term of the order of
$\sqrt{(1/n) \lambdamax\log\log d}$ for all
distributions satisfying the fourth-moment assumption given above.
The additional term is clearly of minor importance. (For example, it
is dominated by the first term whenever $\Tr(\Sigma) > \lambdamax\log\log d$.)
With the estimator
we construct, this term is inevitable.
On the other hand, the inequality of the theorem only holds for sample
sizes that are at least a constant times $d\log d$. 
This feature is not desirable for truly high-dimensional problems, especially taking into account 
that Minsker's bound is ``dimension-free''.

 The fourth-moment assumption can be interpreted as a boundedness assumption of the kurtosis of $(X-\mu)^Tv$. The same assumption has be used in Catoni \cite{catoni2016pac} and Giulini \cite{giulini2015robust} for the robust estimation of the Gram matrix.  The fourth-moment assumption may be weakened to an analogous ``$(2+\epsilon)$-th moment assumption'' that we do not detail for the clarity of the exposition.

We prove the theorem by constructing an estimator in several steps. First
we construct an estimator that performs well for ``spherical''
distributions (i.e., for distributions whose covariance matrix has a trace
comparable to  $d\lambdamax$). This estimator is described in Section \ref{sec:spherical}.
In the second step,
we decompose the space--in a data-dependent way--into the orthogonal
sum of $O(\log d)$ subspaces such that all but one subspaces are such that the projection
of $X$ to the subspace has a spherical distribution. The last subspace
is such that the projection has a covariance matrix with a small trace.
In each subspace
we apply the first estimator and combine them to obtain the final
estimator  $\wh{\mu}_n^{(\delta)}$. 
The proof below provides an explicit value of the constant $C$, though
no attempt has been made to optimize its value.

The constructed estimator is computationally so demanding that even
for moderate values of $d$ it is hopeless to compute it in reasonable
time. 
In this sense, Theorem \ref{thm:main} should be regarded as an
existence result.
It is an interesting an important challenge to construct
estimators with similar statistical performance that can be computed
in polynomial time (as a function of $n$ and $d$). Note that the
estimator of Minsker cited above may be computed by solving a
convex optimization problem, making it computationally feasible, see 
also Hsu and Sabato \cite{HsSa16} for further computational considerations.

\section{An estimator for spherical distributions}
\label {sec:spherical}

In this section we construct an estimator that works well whenever the
distribution of $X$ is sufficiently spherical in the sense that a
positive fraction of the eigenvalues of the covariance matrix is of
the same order as $\lambdamax$. More precisely, for $c \ge 1$, we call a distribution \emph{$c$-spherical} if 
$ d\lambdamax \le c\Tr(\Sigma)$. 

For each $\delta\in (0,1)$ and unit vector $w\in S^{d-1}$ 
(where $S^{d-1}=\{x\in \R^d: \|x\|=1\}$), we may define 
$m_n^{(\delta)}(w)$ as the median-of-means estimate 
(as defined in Section \ref{sec:1d}) of 
$w^T\mu =\EXP w^TX$ based on the i.i.d.\ sample 
$w^TX_1,\ldots,w^TX_n$.

Let $N_{1/2}\subset S^{d-1}$ be a minimal $1/2$-cover, that is, a set of
smallest cardinality that has the property that for all $u\in S^{d-1}$
there exists $w\in N_{1/2}$ with $\|u-x\| \le 1/2$. It is well known 
(see, e.g., \cite[Lemma 13.1.1]{mat02}) that $|N_{1/2}|\le 8^d$.

Noting that $\var(w^T X)\le \lambdamax$,
by (\ref{eq:mom}) and the union bound, we have that, with probability at least $1-\delta$,
\[
 \sup_{w\in N_{1/2}}  \left| m_n^{(\delta/8^d)}(w) - w^T\mu \right| 
\le 2e\sqrt{2  \lambdamax \frac{\ln(e8^d/\delta)}{n}}~.
\]
In other words, if, for $\lambda >0$, we define the empirical polytope
\[
   P_{\delta,\lambda}= \left\{ x\in \R^d: \sup_{w\in N_{1/2}}  \left| m_n^{(\delta/8^d)}(w) - w^Tx \right| 
\le 2e\sqrt{2  \lambda \frac{\ln(e8^d/\delta)}{n}}\right\}~,
\]
then with probability at least $1-\delta$, $\mu \in
P_{\delta,\lambdamax}$. In particular, on this event,
$P_{\delta,\lambdamax}$ is nonempty. Suppose that an upper bound of
the largest eigenvalue of the covariance matrix
$\lambda \ge \lambdamax$ is available. Then we may define the mean
estimator
\[
   \wh{\mu}_{n,\lambda}^{(\delta)}= \left\{ \begin{array}{ll}
                 \text{any element $y\in P_{\delta,\lambda}$} & \text{if
                      $P_{\delta,\lambda}\neq \emptyset$} \\ 
                        0 & \text{otherwise} \end{array} \right. ~.
\]
Now suppose that $\mu \in P_{\delta,\lambda}$ and let $y\in
P_{\delta,\lambda}$ be arbitrary.
Define
$u=(y-\mu)/\|y-\mu\| \in S^{d-1}$, and let $w\in N_{1/2}$ be such that
$\|w-u\|\le 1/2$. (Such a $w$ exists by definition of $N_{1/2}$.) Then
\[
   \|y-\mu\|= u^T(y-\mu)= (u-w)^T(y-\mu) + w^T(y-\mu) \le (1/2)
   \|y-\mu\| + 4e\sqrt{2  \lambda \frac{\ln(e8^d/\delta)}{n}}~,
\]
where we used Cauchy-Schwarz and the fact that $y,\mu\in P_{\delta,\lambda}$.
Rearranging, we obtain that, on the event that $\mu \in P_{\delta,\lambda}$,
\[
\left\|\wh{\mu}_{n,\lambda}^{(\delta)}-\mu\right\| \le  8e\sqrt{2  \lambda
  \frac{d\ln 8+ \ln(e/\delta)}{n}}~,
\]
provided that $\lambda \ge \lambdamax$. Summarizing, we have proved
the following.

\begin{proposition}
\label {prop:spherical}
Let $\lambda>0$ and $\delta \in (0,1)$. 
For any distribution with mean $\mu$ and covariance matrix $\Sigma$
such that $\lambdamax= \|\Sigma\| \le \lambda$, the estimator
$\wh{\mu}_{n,\lambda}^{(\delta)}$ defined above satisfies, with
probability at least $1-\delta$,
\[
\left\|\wh{\mu}_{n,\lambda}^{(\delta)}-\mu\right\| \le  8e\sqrt{2  \lambda
  \frac{d\ln 8+ \ln(e/\delta)}{n}}~.
\]
In particular, if the distribution is $c$-spherical and $\lambda \le 2\lambdamax$, then
\[
\left\|\wh{\mu}_{n,\lambda}^{(\delta)}-\mu\right\| \le  16e\sqrt{  
  \frac{c\Tr(\Sigma)\ln 8+ \lambdamax\ln(e/\delta)}{n}}~.
\]
\end{proposition}

The bound we obtained has the same sub-Gaussian form as
(\ref{eq:subgauss}),
up to a multiplicative constant, whenever the distribution is
$c$-spherical.
To make the estimator fully data-dependent, we need to find an estimate
$\wh{\lambda}$ that falls in the interval $[\lambdamax,2\lambdamax]$,
with high probability. This may be achieved by splitting the sample in 
two parts of equal size (assuming $n$ is even), estimating $\lambdamax$
using samples from one part and computing the mean estimate
defined above using the other part. In the next section we describe
such a method as a part of a more general procedure.

\section{Empirical eigendecomposition}

In the previous section we presented a mean estimate that works
well for ``spherical'' distributions. We will use this estimator as a building
block in the construction of an estimator that has the desirable 
performance guarantee for distributions with any covariance matrix. 
In addition to finite covariances, we assume that there
exists a constant $K>0$ such that, for all $v\in \Rd$ with $\|v\|=1$,
\begin{equation}
\label{eq:4thmoment}
   \EXP\left[ \left((X-\mu)^Tv\right)^4\right] \le K (v^T\Sigma v)^2~.
\end{equation}
In this section we assume that
$n\ge 2(400e)^2K\log_{3/2}d\left(d\log 25 + \log
  (2\log_{3/2}d) + \log (1/\delta)\right)$.

The basic idea is the following. We split the data into two equal
halves. We use the first half in order to decompose the space into
the sum of orthogonal subspaces such that the projection of $X$ into 
each subspace is $4$-spherical. Then we may estimate the projected means
by the estimator of the previous section. 

Next we describe how we obtain an orthogonal decomposition of the
space  based on $n$ i.i.d.\ observations $X_1,\ldots,X_n$. 

Let $s=\lceil \log_{3/2} d^2 \rceil$ and $m= \lfloor n/s \rfloor$.
Divide the sample into $s$ blocks, each of size at least $m$.
In what follows, we describe a way of sequentially decomposing $\Rd$ into the orthogonal sum of 
$s+1$ subspaces $\Rd = V_1 \oplus \cdots \oplus V_{s+1}$. First we construct
$V_1$
using the first block $X_1,\ldots,X_m$ of observations. Then we use
the second block to build $V_2$, and so on, for $s$ blocks. 
The key properties we need
are
that (a) the random vector $X$, projected to any of these subspaces
has a $4$-spherical distribution; (b) the largest eigenvalue of the
covariance matrix of $X$, projected on $V_i$ is at most $\lambdamax (2/3)^{i-1}$.

To this end,
just like in the previous section, let $N_{\gamma}\subset S^{d-1}$ be a 
minimal $\gamma$-cover of the unit sphere $S^{d-1}$ for a sufficiently
small constant
$\gamma\in (0,1)$.  
The value $\gamma=1/100$ is sufficient for our purposes and in the
sequel we assume this value.
Note that $|N_{\gamma}|\le (4/\gamma)^d$ (see \cite[Lemma 13.1.1]{mat02} for a proof of this fact).

Initially, we use the first block $X_1,\ldots,X_m$. We may assume 
that $m$ is even. 
Using these
observations, for each $u\in N_{\gamma}$, we compute an 
estimate $V_m^{(\delta)}(u)$ of 
$u^T\Sigma u=\EXP (u^T(X-\mu))^2=(1/2)\EXP (u^T(X-X'))^2$, where
$X'$ is an i.i.d.\ copy of $X$. We may construct the estimate by
forming $m/2$ i.i.d.\ random variables $(1/2)(u^T(X_1-X_{m/2+1}))^2,\ldots,(1/2)(u^T(X_{m/2}-X_m))^2$
and estimate their mean by the median-of-means estimate $V_m^{(\delta)}(u)$ with parameter $\delta/(s(4/\gamma)^d)$.
Then (\ref{eq:mom}), together with assumption (\ref{eq:4thmoment})
implies that, with probability at least $1-\delta/s$,
\[
    \sup_{u\in N_{\gamma}} \frac{\left|u^T\Sigma u - V_m^{(\delta)}(u)
      \right|}{u^T\Sigma u} \le 4e\sqrt{\frac{K \log(s (4/\gamma)^d/\delta)}{m}}
    \defeq \epsilon_m~.
\]
Our assumptions on the sample size guarantee that $\epsilon_m<1/100$.
The event that the inequality above holds is denoted by $E_1$ so that
$\PROB\{E_1\} \ge 1-\delta/s$.

Let $\M_{\delta,m}$ be the set of all symmetric positive semidefinite $d\times
d$ matrices $M$ satisfying
\[
    \sup_{u\in N_{\gamma}} \frac{\left|u^TM u - V_m^{(\delta)}(u)
      \right|}{u^T\Sigma u} \le \epsilon_m~.
\]
By the argument above, $\Sigma \in \M_{\delta,m}$ on the event
$E_1$. In particular, on $E_1$,
$\M_{\delta,m}$ in non-empty. Define the estimated covariance matrix
\[
   \wh{\Sigma}_{m}^{(\delta)}= \left\{ \begin{array}{ll}
                 \text{any element of $\M_{\delta,m}$} & \text{if
                      $\M_{\delta,m}\neq \emptyset$} \\ 
                        0 & \text{otherwise} \end{array} \right. 
\]
Since on $E_1$ both $\wh{\Sigma}_{m}^{(\delta)}$ and $\Sigma$ are in
$\M_{\delta,m}$, on this event, we have
\begin{equation}
\label {eq:quadforms}
  \left(u^T\Sigma u\right) \frac{1-\epsilon_m}{1+\epsilon_m} \le u^T  \wh{\Sigma}_{m}^{(\delta)}u
\le   \left(u^T\Sigma u \right) \frac{1+\epsilon_m}{1-\epsilon_m}
\qquad \text{for all $u\in N_{\gamma}$}.
\end{equation}

Now compute the spectral decomposition 
\[
   \wh{\Sigma}_{m}^{(\delta)} = \sum_{i=1}^d \wh{\lambda}_i \wh{v}_i \wh{v}_i^T~,
\]
where $\wh{\lambda}_1\ge \cdots \ge \wh{\lambda}_d\ge 0$ are the
eigenvalues and $\wh{v}_1,\ldots,\wh{v}_d$ the corresponding
orthogonal eigenvectors.

Let $u\in S^{d-1}$ be arbitrary and let $v$ be 
a point in $N_{\gamma}$ with smallest distance to $u$. Then
\begin{eqnarray}
  u^T\wh{\Sigma}_{m}^{(\delta)} u & = &
v^T \wh{\Sigma}_{m}^{(\delta)} v + 2(u-v)^T \wh{\Sigma}_{m}^{(\delta)}
v + (u-v)^T\wh{\Sigma}_{m}^{(\delta)} (u-v)  
\nonumber\\
& \le &
v^T \wh{\Sigma}_{m}^{(\delta)} v + \wh{\lambda}_1 (2\gamma + \gamma^2) 
\label {eq:uv1}
\\
& & \text{(by Cauchy-Schwarz and using the fact that $\|u-v\|\le
  \gamma$)}  \nonumber \\
& \le & (v^T \Sigma v) \frac{1+\epsilon_m}{1-\epsilon_m}+ 3
\gamma\wh{\lambda}_1  \nonumber \\
& & \text{(by (\ref{eq:quadforms}))} \nonumber \\
& \le & \frac{1+\epsilon_m}{1-\epsilon_m} \lambdamax + 3\gamma\wh{\lambda}_1~.
\nonumber
\end{eqnarray}
In particular, on $E_1$
we have $\wh{\lambda}_1 \le \beta\lambdamax$ where
$\beta=\frac{1+\epsilon_m}{1-\epsilon_m}/(1-3\gamma) <1.1$.

By a similar argument, we have that for any $u\in S^{d-1}$, 
if $v$ is the point in $N_{\gamma}$ with smallest distance to $u$, then on $E_1$,
\[
  u^T\Sigma u \le 
(v^T \wh{\Sigma}_{m}^{(\delta)} v) \frac{1+\epsilon_m}{1-\epsilon_m} +
3\gamma\lambdamax~
  \le  \frac{1+\epsilon_m}{1-\epsilon_m}\wh{\lambda}_1 +
3\gamma \lambdamax~.
\]
In particular, $\lambdamax \le \beta \wh{\lambda}_1 \le (4/3)\wh{\lambda}_1$.
Similarly,
\begin{eqnarray}
\nonumber
 u^T\Sigma u & \ge &
(v^T \wh{\Sigma}_{m}^{(\delta)} v) \frac{1-\epsilon_m}{1+\epsilon_m} - 3\gamma\wh{\lambda}_1  
\\
\nonumber
& \ge &
\left(u^T \wh{\Sigma}_{m}^{(\delta)} u -3\gamma\wh{\lambda}_1\right) \frac{1-\epsilon_m}{1+\epsilon_m} - 3\gamma\wh{\lambda}_1~.
\\
\label{eq:l1}
& \ge &
\left(u^T \wh{\Sigma}_{m}^{(\delta)} u\right) \frac{1-\epsilon_m}{1+\epsilon_m} - 6\gamma\wh{\lambda}_1~.
\end{eqnarray}

Let $\wh{d}_1$ be number of eigenvalues
$\wh{\lambda}_i$ that are at least $\wh{\lambda}_1/2$ and let $V_1$  
be the subspace of $\R^d$ spanned by
$\wh{v}_1,\ldots,\wh{v}_{\wh{d}_1}$. Denote by $\Pi_1(X)$ the orthogonal projection of
the random variable $X$ (independent of the $X_i$ used to build $V_1$)
onto $V_1$. Then for any $u\in  V_1\cap S^{d-1}$,
 on the event $E_1$, by (\ref{eq:l1}), 
\[
u^T \Sigma u \ge \wh{\lambda}_1\frac{1}{2}\left(\frac{1-\epsilon_m}{1+\epsilon_m} -
12\gamma\right) \ge \frac{\wh{\lambda}_1}{3}
\]
and therefore
\[
   \EXP u^T (\Pi_1(X)-\EXP \Pi_1(X)) (\Pi_1(X)-\EXP \Pi_1(X))^T u = u^T \Sigma u \in
   \left(\frac{\wh{\lambda}_1}{3}, \frac{4\wh{\lambda}_1}{3} \right)~.
\]
In particular, the ratio of the largest and smallest eigenvalues of the
covariance matrix of $\Pi_1(X)$ is at most $4$ and therefore the
distribution of $\Pi_1(X)$ is $4$-spherical. 

On the other hand, on the event $E_1$, for any unit vector $u\in
V_1^{\bot}\cap S^{d-1}$
in the orthogonal complement of $V_1$, we have $u^T \Sigma u \le
2\lambdamax/3$.
To see this, note that
$u^T \wh{\Sigma}_{m}^{(\delta)} u\le \wh{\lambda}_1/2$ and therefore, 
denoting by $v$ the point in $N_{\gamma}$ closest to $u$, 
\begin{eqnarray*}
u^T \Sigma u
& = &
u^T \wh{\Sigma}_{m}^{(\delta)} u +
v^T \left(\Sigma-\wh{\Sigma}_{m}^{(\delta)} \right) v
+ \left(v^T \wh{\Sigma}_{m}^{(\delta)} v- 
 u^T \wh{\Sigma}_{m}^{(\delta)} u \right) +
\left(u^T \Sigma u - v^T \Sigma v \right)  \\
& \le &
\frac{\wh{\lambda}_1}{2}
+ 2\epsilon_m\lambdamax
+ 3\gamma \wh{\lambda}_1  + 3\gamma \lambdamax
\\
& & \text{(by (\ref{eq:quadforms}), (\ref{eq:uv1}), and a similar
  argument for the last term)}  \\
& \le & 
\lambdamax \left(  \beta  \left(\frac{1}{2} +3\gamma\right) +
  2\epsilon_m+3\gamma \right) \le \frac{2\lambdamax}{3}~.
\end{eqnarray*}
In other words, the  largest eigenvalue of the covariance matrix of
$\Pi_1^{\bot}(X)$ (the projection of $X$ to the subspace $V_1^{\bot}$) is at most
$(2/3)\lambdamax$.

In the next step we construct the subspace $V_2 \subset V_1^{\bot}$.
To this end, we proceed exactly as in the first step but now we
replace $\R^d$ by $V_1^{\bot}$ and the sample $X_1,\ldots,X_m$ on the
first block by the variables $\Pi_1^{\bot}(X_{m+1}),\ldots,\Pi_1^{\bot}(X_{2m})\in
V_1^{\bot}$. (Recall that $\Pi_1^{\bot}(X_i)$ is the projection of $X_i$ to
the subspace $V_1^{\bot}$). Just like in the first step, with
probability at least $1-\delta/s$ we obtain a (possibly empty)
subspace $V_2$, orthogonal to $V_1$ such that $\Pi_2(X)$, the projection of $X$
on $V_2$, has a $4$-spherical distribution and largest eigenvalue of the covariance matrix of
$\Pi_2^{\bot}(X)$ (the projection of $X$ to the subspace
$(V_1\oplus V_2)^{\bot}$) is at most
$(2/3)^2\lambdamax$.

We repeat the procedure $s$ times and use a union bound the $s$ events. We obtain, with probability at
least $1-\delta$, a sequence of subspaces $V_1,\ldots,V_s$, with the
following properties:

\begin{enumerate}
\item[(i)] $V_1,\ldots,V_s$ are orthogonal subspaces.
\item[(ii)]
For each $i=1,\ldots,s$, 
$\Pi_i(X)$, the projection of $X$
on $V_i$, has a $4$-spherical distribution.
\item[(iii)]
The largest eigenvalue of the covariance matrix of
$\Pi_i(X)$ is at most
$\lambda_1^{(i)}\le (2/3)^{i-1}\lambdamax$.
\item[(iv)]
The largest eigenvalue $\wh{\lambda}_1^{(i)}$ of the estimated
covariance matrix of $\Pi_i(X)$ satisfies
\[
 (3/4) \lambda_1^{(i)}   \le  \wh{\lambda}_1^{(i)} \le 1.1 \lambda_1^{(i)}~.
\]
\end{enumerate}

Note that it may happen for some $T<s$, we have $\Rd=V_1 \oplus \cdots
\oplus V_T$. In that case we define $V_{T+1}= \cdots =V_s = \emptyset$.

\section{Putting it all together}

In this section we construct our final multivariate mean estimator and
prove Theorem \ref{thm:main}. To simplify notation, we assume that the
sample size is $2n$. This only effects the value of the universal
constant $C$ in the statement of the theorem.

The data is split into two equal halves $(X_1,\ldots,X_n)$ and
$(X_{n+1},\ldots,X_{2n})$. 
The second half is used to construct the orthogonal spaces
$V_1,\ldots,V_s$ as described in the previous section. Let
$\wh{d}_1,\ldots,\wh{d}_s$ denote the dimension of these subspaces.
Recall that, with probability at least $1-\delta$, the construction is 
successful in the sense that the subspaces satisfy properties  (i)--(iv)
described at the end of the previous section. Denote this event by $E$.
In the rest of the argument we condition on $(X_{n+1},\ldots,X_{2n})$
and assume that $E$ occurs. All probabilities below are conditional.

 If 
$\sum_{i=1}^s \wh{d}_i < d$ (i.e., $V_1\oplus\cdots \oplus V_s \neq
\Rd$),
then we define $V_{s+1} = (V_1\oplus\cdots \oplus V_s)^{\bot}$ and denote
by $\wh{d}_{s+1} = d- \sum_{i=1}^s \wh{d}_i$ the dimension of
$V_{s+1}$.
Let $\Pi_1,\ldots,\Pi_{s+1}$ denote the projection operators on the
subspaces $V_1,\ldots,V_{s+1}$, respectively.
For each $i=1,\ldots,s+1$, we use the vectors
$\Pi_i(X_1),\ldots,\Pi_i(X_n)$ to compute an estimator of the mean
$\EXP \left[\Pi_i(X)|(X_{n+1},\ldots,X_{2n})\right] = \Pi_i(\mu)$. 

For $i=1,\ldots,s$, we use the estimator defined in Section \ref{sec:spherical}.
In particular, within the $\wh{d}_i$-dimensional space $V_i$, we
compute $\ol{\mu}_i =
\wh{\mu}_{n,(4/3)\wh{\lambda}_i}^{(\delta/(s+1))}$. 
Note that since $\wh{\lambda}_i$ comes from an empirical estimation of $\Sigma$ restricted to an empirical subspace $V_i$, $\ol{\mu}_i$ is an estimator constructed on the sample $X_1,\dots,X_n$.
Then, by Proposition \ref {prop:spherical},
with probability $1-\delta/(s+1)$,
\[
\left\|\ol{\mu}_i-\Pi_i(\mu)\right\|^2 \le  (8e)^2  
  \frac{ (8/3) \wh{\lambda}_1^{(i)} \left(\wh{d}_i \ln 8+ \ln(e (2\log_{3/2}d+1) /\delta)\right)}{n}~.
\]
In the last subspace $V_{s+1}$, we may use Minsker's estimator, based
on
$\Pi_{s+1}(X_1),\ldots,\Pi_{s+1}(X_n)$ to compute an estimator
$\ol{\mu}_{s+1}= \wt{\mu}_n^{(\delta/(s+1))}$  of
$\Pi_{s+1}(\mu)$. Since the largest eigenvalue of the covariance
matrix of $\Pi_{s+1}(X)$ is at most $\lambdamax/d^2$, using
(\ref{eq:minsker}), we obtain that, with probability $1-\delta/(s+1)$,
\[
 \|\ol{\mu}_{s+1}-\Pi_{s+1}(\mu)\|^2 \le C\frac{\lambdamax \log((2\log_{3/2}d+1)/\delta)}{n}~.
\]
Our final estimator is $\wh{\mu}_n^{(\delta)} = \sum_{i=1}^{s+1}
\ol{\mu}_{s+1}$. By the union bound, we have that, with
probability at least $1-\delta$,
\begin{eqnarray*}
\left\|  \wh{\mu}_n^{(\delta)} -\mu \right\|^2
& = & \sum_{i=1}^{s+1} \left\|\ol{\mu}_i-\Pi_i(\mu)\right\|^2 \\
& \le &  (8e)^2  
  \frac{ (8/3) \ln 8 }{n} \sum_{i=1}^s  \wh{\lambda}_1^{(i)}\wh{d}_i+ 
(8e)^2(8/3) \frac{\ln(e    (2\log_{3/2}d+1) /\delta)}{n}  \sum_{i=1}^s
 \wh{\lambda}_1^{(i)}  \\
& & \qquad + C\frac{\lambdamax \log((2\log_{3/2}d+1)/\delta)}{n}
\end{eqnarray*}
First notice that, by properties (iii) and (iv) at the end of the
previous section, 
\[
 \sum_{i=1}^s  \wh{\lambda}_1^{(i)}  \le 1.1  \sum_{i=1}^s
 \lambda_1^{(i)} 
\le 1.1  \lambdamax  \sum_{i=1}^s  (2/3)^{i-1} \le 3.3 \lambdamax~.
\]
On the other hand, since 
\[
\Tr(\Sigma) = \EXP \|X-\mu\|^2 = \sum_{i=1}^{s+1} \EXP\| \Pi_i(X)-\Pi_i(\mu)\|^2
\]
and for $i\le s$ each $\Pi_i(X)$ has a
$4$-spherical distribution, we have that
\[
  \sum_{i=1}^s  \wh{\lambda}_1^{(i)} \wh{d}_i \le 1.1   \sum_{i=1}^s
  \lambda_1^{(i)} \wh{d}_i 
  \le 4.4 \Tr(\Sigma)~. 
\]
This concludes the proof of Theorem \ref{thm:main}.

\bibliographystyle{plain}
\bibliography{meanest}

\end{document}